\documentclass[10pt]{article}

\usepackage{amscd,amsmath, amssymb, fancyhdr}

\numberwithin{equation}{section}

\newcommand{\version}{version 1.4,\ \  30 April 2008}

\def\eqref#1{(\ref{#1})}
\newcommand{\goth}{\mathfrak}
\newcommand{\g}{{\mathfrak g}}
\newcommand{\arrow}{{\:\longrightarrow\:}}
\newcommand{\Z}{{\Bbb Z}}
\newcommand{\C}{{\Bbb C}}
\newcommand{\R}{{\Bbb R}}

\renewcommand{\H}{{\Bbb H}}
\newcommand{\6}{\partial}
\def\1{\sqrt{-1}\:}

\renewcommand{\tilde}{\widetilde}
\renewcommand{\bar}{\overline}
\renewcommand{\phi}{\varphi}
\renewcommand{\epsilon}{\varepsilon}

\renewcommand{\leq}{\leqslant}

\newcommand{\End}{\operatorname{End}}

\newcommand{\Id}{\operatorname{Id}}

\newcommand{\Hol}{\operatorname{Hol}}

\newcommand{\Alt}{\operatorname{Alt}}

\newcounter{Mycounter}[section]
\newcounter{lemma}[section]
\renewcommand{\thelemma}{{Lemma \thesection.\arabic{lemma}}}
\newcommand{\lemma}{%
  \setcounter{lemma}{\value{Mycounter}}
  \refstepcounter{lemma}
  \stepcounter{Mycounter}
  {\noindent \bf \thelemma:\ }}

\newcounter{claim}[section]
\newcounter{sublemma}[section]
\newcounter{corollary}[section]
\renewcommand{\thecorollary}{{Corollary \thesection.\arabic{corollary}}}
\newcommand{\corollary}{%
  \setcounter{corollary}{\value{Mycounter}}
  \refstepcounter{corollary}
  \stepcounter{Mycounter}
  {\noindent \bf \thecorollary:\ }}

\newcounter{theorem}[section]
\renewcommand{\thetheorem}{{Theorem \thesection.\arabic{theorem}}}
\newcommand{\theorem}{%
  \setcounter{theorem}{\value{Mycounter}}
  \refstepcounter{theorem}
  \stepcounter{Mycounter}
  {\noindent \bf \thetheorem:\ }}

\newcounter{conjecture}[section]
\newcounter{proposition}[section]
\renewcommand{\theproposition}
    {{Proposition \thesection.\arabic{proposition}}}
\newcommand{\proposition}{%
  \setcounter{proposition}{\value{Mycounter}}
  \refstepcounter{proposition}
  \stepcounter{Mycounter}
  {\noindent \bf \theproposition:\ }}

\newcounter{definition}[section]
\renewcommand{\thedefinition}
    {{Definition~\thesection.\arabic{definition}}}
\newcommand{\definition}{%
  \setcounter{definition}{\value{Mycounter}}
  \refstepcounter{definition}
  \stepcounter{Mycounter}
  {\noindent \bf \thedefinition:\ }}

\newcounter{example}[section]
\newcounter{remark}[section]
\renewcommand{\theremark}{{Remark \thesection.\arabic{remark}}}
\newcommand{\remark}{%
  \setcounter{remark}{\value{Mycounter}}
  \refstepcounter{remark}
  \stepcounter{Mycounter}
  {\noindent \bf \theremark:\ }}

\newcounter{problem}[section]
\newcounter{question}[section]
\makeatletter

\@addtoreset{equation}{section} \@addtoreset{footnote}{section}
\makeatother

\def\blacksquare{\hbox{\vrule width 5pt height 5pt depth 0pt}}
\def\endproof{\blacksquare}

\begin{document}
\begin{center}
{\LARGE\bf
Canonical bundles of  complex \\[2mm]nilmanifolds, with applications
to\\[2mm] hypercomplex geometry\\[5mm]
}

Mar\'ia L. Barberis, Isabel G. Dotti and Misha
Verbitsky\footnote{Misha Verbitsky is partially 
supported by a CRDF grant RM1-2354-MO02.}

\end{center}

{\small \hspace{0.15\linewidth}
\begin{minipage}[t]{0.7\linewidth}
{\bf Abstract} \\
A nilmanifold is a quotient of a nilpotent group $G$ by a co-compact
discrete subgroup. A complex nilmanifold is one which is equipped
with a $G$-invariant complex structure. We prove that a complex
nilmanifold has trivial canonical bundle. This is used to study
hypercomplex nilmanifolds (nilmanifolds with a triple of
$G$-invariant complex structures which satisfy quaternionic
relations). We prove that a hypercomplex nilmanifold admits an HKT
(hyperk\"ahler with torsion) metric if and only if   the underlying
hypercomplex structure is abelian. Moreover, any $G$-invariant
HKT-metric on a nilmanifold is balanced with respect to all
associated complex structures.

\end{minipage}
}

\tableofcontents


\section{Introduction}


\subsection{Canonical bundle of complex nilmanifolds}

Let $G$ be a connected, simply connected nilpotent Lie group, and
$\Gamma\subset G$ a discrete, co-compact subgroup. The quotient
manifold $\Gamma \backslash G$ is called {\bf a nilmanifold}.
Clearly, $G$ acts on $\Gamma \backslash G$ transitively (from the
right). Nilmanifolds are often defined as compact manifolds with a
transitive action of a nilpotent Lie group (see e.g.
\cite{_Hasegawa_}). In this case, the above definition becomes a
theorem, proven by Mal'cev, who introduced this notion in 1949, in
the influental paper \cite{_Malcev_}.

If $N= \Gamma \backslash G$ is equipped with a complex structure
$\mathcal I$ induced by a left-invariant complex structure  on $G$,
$(N,{\mathcal I})$ is called {\bf a complex nilmanifold}.

%
%

It is important to note that $G$ is not necessarily a complex Lie
group. Indeed, for $G$ to be a complex Lie group, both left and
right translations on $G$ must be holomorphic. In many examples of
complex nilmanifolds, even the simplest ones (such as a Kodaira
surface), this condition is not satisfied.

Geometry of nilmanifolds is an important subject, much studied since
Mal'cev and Mostow started this work. The complex nilmanifolds are
much less understood. There are many papers dealing with special
cases of nilpotent Lie groups and particular discrete co-compact
subgroups, where  the complex invariants of the corresponding
nilmanifolds (deformation space, Fr\"olicher spectral sequence, and
so on) have been computed. However, general results about complex
nilmanifolds remain scarce. Indeed, nilmanifolds serve mostly as a
rich source of counterexamples to all kinds of general conjectures
in geometry and topology (\cite{_Abbena_}, \cite{_CFG:_symple_},
\cite{_CFL:_comple_}).

In the present paper, we prove that the canonical bundle of any
complex nilmanifold is holomorphically trivial
(\ref{_canoni_trivial_Theorem_}). This condition is quite strong.
For instance, any compact complex surface with trivial canonical
bundle is isomorphic to a K3 surface, a torus, or a Kodaira surface;
the first two are K\"ahler, and the latter is a nilmanifold.

K\"ahler manifolds with trivial canonical bundle play an important
role in mathematics and physics, due to the nice behavior of their
deformation spaces, guaranteed by Bogomolov-Tian-Todorov theorem.
This theorem states that the deformations of a compact K\"ahler
manifold with trivial canonical bundle are non-obstructed, and the
deformation space is smooth. 

For nilmanifolds obtained as quotients of complex
nilpotent groups (``complex parallelisable
nilmanifolds''), this is known to be false,
as S. Rollenske proved (\cite{_Rollenske:Kuranishi_}).
For other classes of nilmanifolds, for instance
hypercomplex nilmanifolds, some version of 
Bogomolov-Tian-Todorov theorem could still be true.
In fact, the key ingredient of the proof
of  Bogomolov-Tian-Todorov theorem, the so-called Tian-Todorov lemma
(\cite{_Barannikov_Kontse_}), remains valid for any complex manifold
with trivial canonical bundle. The rest of the proof, sadly, does not
work, because it requires the degeneration of a Fr\"olicher spectral
sequence, and (as shown in \cite{_CFG:_Frolicher_} and
\cite{_Rollenske:Frolicher_}), this sequence can be arbitrarily
non-degenerate. Still, the vast empirical evidence (see e.g.
\cite{_MPP:Deformations_}, \cite{_Poon:Kodaira_defo_}) shows that
some analogue of Tian-Todorov formalism could exist on some (or all)
nilmanifolds.

For general compact non-K\"ahler manifolds with trivial canonical
bundle, an analogue of Bogomolov-Tian-Todorov theorem is known to
be false. In \cite{_Ghys_}, it was shown that the deformation
space of a locally homogeneous manifold $ SL(2,\C)/\Gamma$ can be
obstructed, for a cocompact and discrete subgroup $\Gamma \subset
SL(2,\C)$. 

One of the first examples
of a complex manifold with obstructed deformations
was constructed by A. Douady, in \cite{_Douady_}.
Douady used an Iwasawa manifold 
which is a quotient $M:=G/\Gamma$, with
$G$ the group of complex upper triangular
$3\times 3$-matrices, and $\Gamma$ the group of
upper triangular matrices with coefficients
in Gaussian integers.

Douady proved that a product $M\times \C P^1$
has obstructed deformation space. In
\cite{_Retakh_}, Douady's construction
was generalized, using Massey operations
on cohomology of $M$.

Another proof of triviality of the canonical bundle 
of a nilmanifold is given in \cite{_CG:generalized_on_nilm_}.

\subsection{Hypercomplex nilmanifolds}

An almost hypercomplex manifold is a smooth manifold $M$ equipped
with three operators $\mathcal{I, J, K}\in \End(TM)$ satisfying the
quaternionic relations $\mathcal{I}\circ \mathcal{J} = -
\mathcal{J}\circ \mathcal{I} =\mathcal{K}$, $\mathcal{I}^2 =
\mathcal{J}^2=\mathcal{K}^2=-\Id_{TM}$. The operators $\mathcal{I}$,
$\mathcal{J}$, $\mathcal{K}$ define almost complex structures on
$M$; if these almost complex structures are integrable, $M$ is
called {\bf hypercomplex}. A hypercomplex manifold is equipped with
a whole 2-dimensional sphere of complex structures.

Hypercomplex manifolds were defined by C.P. Boyer (\cite{_Boyer_}),
who gave a classification of compact hypercomplex manifolds for
$\dim_{\Bbb H} M =1$. Many interesting examples of hypercomplex
manifolds were found in the 90-ies, see e.g. \cite{_Joyce_},
\cite{_Pedersen_Poon:inhomo_}, \cite{_Barberis_Dotti_}.
Independently (and earlier) some of these constructions were
obtained by string physicists; see e.g. \cite{_SSTvP_}.

As Obata has shown (\cite{_Obata_}), a hypercomplex manifold
admits a (necessarily unique) torsion-free connection, preserving
$\mathcal{I,J,K}$. The converse is also true: if an almost
hypercomplex manifold admits a torsion-free connection preserving
the quaternionic action, it is hypercomplex. This implies that a
hypercomplex structure on a manifold can be defined as a
torsion-free connection with holonomy in $GL(n, {\Bbb H})$.

Connections with restricted holonomy is one of the central notions
in Riemannian geometry, due to Berger's classification of
irreducible holonomy of Riemannian manifolds. However, a similar
classification exists for a general torsion-free connection
(\cite{_Merkulov_Sch:long_}). In the Merkulov-Schwachh\"ofer list,
only three subroups of $GL(n, {\Bbb H})$ occur. In addition to the
compact group $Sp(n)$ (which defines hyperk\"ahler geometry), also
$GL(n, {\Bbb H})$ and its commutator $SL(n, {\Bbb H})$ appear,
corresponding to hypercomplex manifolds and hypercomplex manifolds
with trivial determinant bundle, respectively. Both of these
geometries are interesting, rich in structure and examples, and
deserve detailed study.

Not much is known about $SL(n, {\Bbb H})$-manifolds. It is easy to
see that $(M,\mathcal{I})$ has holomorphically trivial canonical
bundle, when $(M,\mathcal{I, J, K})$ is a hypercomplex manifold with
holonomy in $SL(n, {\Bbb H})$ (\cite{_Verbitsky:canoni_}). For a
hypercomplex $SL(n, {\Bbb H})$-manifold admitting a special kind of
quaternionic Hermitian metric called HKT metric, a version of Hodge
theory was constructed (\cite{_Verbitsky:HKT_}). Using this result,
it was shown that a compact hypercomplex manifold with trivial
canonical bundle has holonomy in $SL(n,{\Bbb H})$, if it admits an
HKT-structure (\cite{_Verbitsky:canoni_})

It is not clear whether the last condition is necessary: for all
known examples of hypercomplex manifolds with trivial canonical
bundle, holonomy lies in $SL(n,{\Bbb H})$.

In the present paper, we prove that holonomy $\Hol(\nabla)$ of a
hypercomplex nilmanifold always lies in $SL(n, {\Bbb H})$
(\ref{_canon_tri_then_SL_Theorem_})

As shown in \cite{AM}, locally $\Hol(\nabla)\subset SL(n, {\Bbb H})$
is equivalent to vanishing of the Ricci curvature of $\nabla$.
However, the vanishing of Ricci curvature is weaker than
$\Hol(\nabla)\subset SL(n, {\Bbb H})$. Consider for example the Hopf
manifold $H= {\Z} \backslash  ({\Bbb H}^n - 0)$. The Obata
connection on $H$ is obviously flat, hence the Ricci curvature
vanishes. However,  $\Hol(\nabla)$ does not lie in $SL(n, {\Bbb H})$.
In fact, the canonical bundle of $H$ is holomorphically non-trivial,
and has no non-zero sections (see Subsection
\ref{_canoni_nilma_Subsection_}).

We give an independent proof of vanishing of Ricci curvature of a
hypercomplex nilmanifold (Section \ref{hypercomplex_holonomy}).

\subsection{Abelian complex structures}

A complex nilmanifold $(N, {\mathcal I})$, with $N=\Gamma\backslash
G$, gives rise to a splitting \[ \g\otimes _{\Bbb R}\Bbb C =
\g^{0,1}\oplus \, \frak g^{1,0},\] where $\frak g^{0,1},\, \frak
g^{1,0}$ are the eigenspaces of the induced complex structure on the
Lie algebra $\frak g$ of $G$.
 By
Newlander-Nirenberg theorem, the almost complex structure ${\mathcal
I}$ is integrable if and only if $ \frak g^{1,0}$ is a complex
subalgebra of $\g\otimes _{\Bbb R}\Bbb C$. $(N, {\mathcal I})$
 is called {\bf abelian}
if the Lie subalgebra $\frak g^{1,0}$ is abelian.

Abelian complex structures were introduced in \cite{_Barberis:PhD_},
and much studied since then (see, for example, \cite{BDM},
\cite{_Barberis_Dotti_}). There are strong restrictions to the
existence of such structures.  In fact, it has been shown by
\cite{P} that the Lie algebra must be two-step solvable. However, a
complete classification is still unknown, though there exist some
partial results (\cite{BD2}). The complex geometry of nilmanifolds
with abelian complex structures is much more accessible than the
general case. In particular, the Dolbeault cohomology of an abelian
nilmanifold can be expressed in terms of the corresponding Lie
algebra cohomology (\cite{_Console_Fino_}, \cite{CFGU}), and the
same is true for the deformation space (\cite{_MPP:Deformations_},
\cite{_Verbitsky:canoni_}, \cite{_Console_Fino_Poon_}).

This notion is specially convenient when applied to hypercomplex
nilmanifolds. If $(N, \mathcal{I, J, K})$ is a hypercomplex
nilmanifold, abelianness of the complex structure $\mathcal{I}$ is
equivalent to the abelianness of $\mathcal{J}$ and $\mathcal{K}$
(\cite{DF3}). Some results on abelian hypercomplex structures can be
found in \cite{_Dotti_Fino:8-dim_}, \cite{B4}.

\subsection{HKT-structures on nilmanifolds}

Let $(M,\mathcal{ I,J,K})$ be a hypercomplex manifold. A
``hyperk\"ahler with torsion'' (HKT) metric on $M$ is a special kind
of a quaternionic Hermitian metric, which became increasingly
important in mathematics and physics during the last seven years.
HKT-metrics were introduced by P. S. Howe and G. Papadopoulos
(\cite{_Howe_Papado_}) and much discussed in the physics and
mathematics literature since then. See \cite{_Gra_Poon_} for a
treatment of HKT-metrics written from a mathematical point of view.
The term ``hyperk\"ahler metric with torsion'' is actually
misleading, because an HKT-metric is not hyperk\"ahler. This is why
we prefer to use the abbreviation ``HKT-manifold''.

 A quaternionic
Hermitian metric is a Riemannian metric which is Hermitian under
$\mathcal{I}$, $\mathcal{J}$ and $\mathcal{K}$. There are three
Hermitian forms associated with such a metric $g$:
\[
\omega_\mathcal{I}= g(\cdot, \mathcal{I}\cdot), \ \
\omega_\mathcal{J}= g(\cdot, \mathcal{J}\cdot), \ \
\omega_\mathcal{K}= g(\cdot, \mathcal{K}\cdot).
\]
When these forms are closed $(M, \mathcal{I, J, K}, g)$ is called a
hyperk\"ahler manifold. In this case, $M$ is also holomorphically
symplectic; indeed, the form $\omega_\mathcal{J} +
\1\omega_\mathcal{K}$ lies in $\Lambda^{2,0}(M,\mathcal{I})$. Being
closed, this $(2,0)$ form is necessarily holomorphic.

The converse is also true: by Calabi-Yau theorem
(\cite{_Besse:Einst_Manifo_}, \cite{_Yau:Calabi-Yau_}), a compact
holomorphically symplectic K\"ahler manifold admits a hyperk\"ahler
metric, which is unique in a given K\"ahler class. In algebraic
geometry, the word ``hyperk\"ahler'' is often used as a synonym to
``holomorphically symplectic''.

The condition $d(\omega_\mathcal{J} + \1\omega_\mathcal{K})=0$ is
equivalent to hyperk\"ahlerianness. A weaker condition
\begin{equation}\label{_HKT_Equation_}
\6(\omega_\mathcal{J} + \1\omega_\mathcal{K})=0
\end{equation}
is often more useful. A quaternionic Hermitian metric $g$ which
satisfies \eqref{_HKT_Equation_} is called HKT (hyperk\"ahler with
torsion). As in the K\"ahler case,  an HKT metric locally has a
potential (see \cite{_Banos_Swann_}).

For abelian hypercomplex nilmanifolds, any left-invariant
quaternionic Hermitian metric is automatically HKT
(\cite{_Dotti_Fino:HKT_}) and  for 2-step nilmanifolds a converse
result was proven in \cite{_Dotti_Fino:HKT_}. Using the triviality
of the canonical bundle and the hypercomplex version of Hodge
theory \cite{_Verbitsky:HKT_}, we generalize the previous result,
showing that any hypercomplex nilmanifold which admits a
left-invariant HKT-metric is in fact abelian
(\ref{_HKT_implies_abelian_}). In \cite{_Gra_Poon_} the question
whether any compact hypercomplex manifold admits an HKT metric was
posed. In particular, a negative answer to this question is given
by a non-abelian hypercomplex nilmanifold, since it has been shown
in \cite{_Fino_Gra_} that existence of any HKT-metric compatible
with a left-invariant hypercomplex structure implies existence of
a left-invariant one. In \S\ref{examples} a family of non-abelian
hypercomplex nilmanifolds is exhibited (see also the nilmanifold
considered in the Remark of \S 4 in \cite{_Dotti_Fino:8-dim_} and
Lemma 3.1 in \cite{_Fino_Gra_}).

We also obtain, as a consequence of \ref{_HKT_implies_abelian_} and
\ref{quat_balanced}, that any invariant HKT-metric on a hypercomplex
nilmanifold is balanced with respect to all underlying complex
structures.


\section{Geometry of complex nilmanifolds}


\subsection{Complex nilmanifolds: basic properties}

\definition
A {\bf nilmanifold} is a quotient $\Gamma \backslash G$ of a
connected simply connected nilpotent Lie group $G$ by a co-compact
discrete subgroup $\Gamma$.

\hfill

By Mal'\v{c}ev theorem (\cite{_Malcev_}), for any simply connected
nilpotent Lie group $G$ with rational structure constants there is a
lattice $\Gamma$ of maximal rank.

\hfill

 Let $G$ be a real Lie group, equipped with a left-invariant
almost complex structure ${\cal I}$, acting on its Lie algebra as
$I:\; {\goth g}\arrow {\goth g}$, $I^2=-\Id$. It is well known that
${\cal I}$ is integrable if and only if the $\1$-eigenspace ${\goth
g}^{0,1}\subset \frak g _{\C} := {\goth
 g}\otimes_\R \C$ is a subalgebra of $\frak g _{\C}$.
In this situation, we shall say that $G$ is equipped with a
left-invariant complex structure. When   $I:\; {\goth g}\arrow
{\goth g}$ satisfies the condition $[Ix,Iy]=[x,y]$ for any $x, y\in
{\goth g}$, ${\cal I}$ is integrable and it is called an {\bf
abelian} complex structure. In this case, it turns out that ${\goth
g}^{1,0}\subset {\goth
 g}\otimes_\R \C$ is a complex abelian subalgebra of ${\goth
 g}\otimes_\R \C$.

\hfill

Let $G$ be a nilpotent Lie group with a left-invariant complex
structure $\mathcal I$. According to  Theorem 1.3 in \cite{_Sal_},
there exist left-invariant $(1,0)$-forms $\omega _1, \dots , \omega
_n$ and smooth $1$-forms $\eta_1^i, \dots , \eta ^i _{i-1}$ on $G$
for $ 2\leq i \leq n$, such that
\begin{equation} \label{Salamon_basis}
 d \omega _i= \sum_{j<i}  \eta _j ^i \wedge \omega _j   .
\end{equation}

\hfill

We prove next an algebraic lemma, which will be useful to prove that
a hypercomplex nilmanifold is Ricci flat (see \ref{ricciflat}). Its
proof makes use of the existence of the above basis of
$(1,0)$-forms.

\hfill

\lemma \label{ric0} Let $\cal J$ be a complex structure on a
nilpotent Lie algebra $\frak g$. Then $$\text{tr}\left( {\cal J}\,
\text{ad}_{X } \right)=0, \text{ for any }X\in \frak g.
$$

\hfill

\noindent {\bf  Proof:} Let \[ \omega _1, \dots , \omega _n \in
\Lambda ^{1,0} \, \frak g \, \] satisfy \eqref{Salamon_basis}, and
consider $\bar{\omega }_1, \dots , \bar{\omega}_n \in \Lambda ^{0,1}
\, \frak g$. If $X_1, \dots, X_n , \bar{X}_1, \dots, \bar{X}_n$ is
the basis of $\frak g_{\C} ^*$ dual to $\, \omega _1, \dots , \omega
_n, \bar{\omega }_1, \dots , \bar{\omega}_n$, then the matrix of
ad$_{X_k}$ relative to this basis takes the form:
\[    \begin{pmatrix} A_k& *\\ 0& B_k \end{pmatrix} ,     \]
where tr$(A_k)=0$ and $B_k$ is strictly lower triangular. In fact,
let $B_k=(b_{il}^k)$. Using \eqref{Salamon_basis} one obtains
\[ {d}\, \bar{\omega}
_i = \sum_{j<i} \bar{\eta}^i_j \wedge  \bar{\omega }_j, \]
 then \begin{eqnarray*} b_{il }^k&=& \bar{\omega} _i \left( [X_k,
\bar{X}_l]\right)=-2 \, {d}\, \bar{\omega}
_i(X_k, \bar{X}_l)\\
& =& -2 \sum _{j<i} \frac12\left( \bar{\eta}^i_j (X_k) \bar{\omega}
_j (\bar{X}_l ) - \bar{\eta}^i_j (\bar{X}_l) \bar{\omega} _j (X_k
)\right)\\
&=& - \sum _{j<i}  \bar{\eta}^i_j (X_k) \bar{\omega} _j (\bar{X}_l
),
\end{eqnarray*}
since $\bar{\omega} _j (X_k )=0$ for any $j, k$. Observe that when
$i\leq l$, $\bar{\omega} _j (\bar{X}_l )=0$ for all $j<i$, therefore
$b_{il}^k=0$ for $i\leq l$, and it turns out that $B_k$ is strictly
lower triangular, as claimed. This implies that tr$(A_k)=0$ since
ad$_{X_k}$ is nilpotent.

On the other hand, the matrix of $\cal J$ relative to $X_1, \dots,
X_n , \bar{X}_1, \dots, \bar{X}_n$ is given by:
\[  \begin{pmatrix} i \, \text{Id}& 0\\ 0& -i\, \text{Id} \end{pmatrix} , \]
therefore, the matrix of ${\cal J}\, \text{ad}_{X _k}$ takes the
following form:
\[ \begin{pmatrix} i A_k& *\\ 0& -i B_k \end{pmatrix} ,\]
and, in particular, it has zero trace. A similar argument, using
that the matrix of ad$_{\bar{X}_k}$  is given by:
\[    \begin{pmatrix} C_k& 0\\ *& D_k \end{pmatrix} ,     \]
with $C_k$ strictly lower triangular and tr$(D_k)=0$, gives that
tr$\left({\cal J}\, \text{ad}_{\bar{X} _k}\right)=0$. Therefore,
tr$\left({\cal J}\, \text{ad}_{{X} }\right)=0$ for any $X\in \frak g
_{\C}$ and the lemma follows.
\endproof

\hfill

Let $(M, {\mathcal J})$ be a complex manifold, $g$ a Hermitian
metric, $\omega =g( \cdot ,\mathcal J \cdot  ) $ the K\"ahler form
and $\theta =d^*\omega \circ \mathcal J$ the Lee form of the
Hermitian manifold $(M, {\mathcal J}, g)$, where $d^*$ is the
adjoint of $d$.

\hfill

\definition A Hermitian metric $g$ on a complex manifold $(M, {\mathcal
J})$ is called {\bf balanced} if $\theta =0$, where $\theta $ is the
associated Lee form.

\hfill

On any Hermitian manifold $(M,\mathcal J,g)$ there exists a unique
connection $\nabla ^B$ satisfying $\nabla ^B g = 0, \; \nabla ^BJ =
0$ and whose torsion tensor $c$ (considered as a $(3,0)$-tensor) $ c
(X,Y,Z) = g (X, T(Y,Z)) $ is totally skew-symmetric. Physicists call
this connection  a $KT$-connection; among mathematicians this
connection is known as the Bismut connection \cite{Bi}. The Lee form
can be expressed (locally) in terms of the torsion tensor $c$ as
follows (see \cite{IP}):
\begin{equation}\label{Lee}
\theta(X)=-\frac 12\sum _{i=1}^{2n} c({\mathcal J}X, E_i , {\mathcal
J} E_i  ),
\end{equation}
for an orthonormal basis $E_1, \dots , E_{2n}$ of (local) vector
fields.

We restrict next to the case of a left invariant Hermitian
structure on a  Lie group. The proof of the next lemma follows by
using the properties of the Bismut connection together with
\eqref{Lee}.

\hfill

\lemma \label{theta} Let $G$ be a  Lie group with an abelian complex
structure ${\mathcal J}$ and $g$ an arbitrary Hermitian
left-invariant metric. Then the Bismut connection $\nabla ^B$ and
the Lee form $\theta$ associated to $(G,{\mathcal J}, g)$ are given
by
\begin{equation}\label{abel_Lee}\begin{split}
g(\nabla ^B_XY,Z)&=-g(X,[Y,Z]),\\  \theta (X)&= \text{tr} \left(
\frac 12 \, J \nabla^B_{JX}-\text{ad}_X\right),
\end{split}\end{equation}
where $X,Y,Z$ are left-invariant vector fields.

\hfill

\subsection{Canonical bundle of a complex nilmanifold}
\label{_canoni_nilma_Subsection_}

\definition
Let $N=\Gamma \backslash G$ be a nilmanifold and assume that $G$ is
equipped with a left-invariant complex structure. This makes $N$
into a complex manifold. In such a situation we say that $N$ is a
{\bf complex nilmanifold}.

\hfill

\definition
A  complex structure ${\cal I}$ on a nilmanifold $N=\Gamma
\backslash G$ is called {\bf abelian } if it is induced from a left
invariant abelian complex structure on $G$.

\hfill

\theorem \label{_canoni_trivial_Theorem_} Let $N= \Gamma \backslash
G$ be a complex nilmanifold, $n=\dim_\C G$. Then $G$ admits a
left-invariant, non-zero, holomorphic section of the canonical
bundle $\Lambda^{n,0}(G)$. In particular, the canonical bundle
$K(N)$ of $N$ is trivial, as a holomorphic line bundle.

\hfill

{\bf Proof:}  Let $\omega _1, \dots , \omega _n$ be the
left-invariant $(1,0)$-forms and  $\eta_1^i, \dots , \eta ^i _{i-1}$
the smooth $1$-forms ($ 2\leq i \leq n$), as in
\eqref{Salamon_basis}. If $\eta = \bigwedge_{i=1} ^n \omega _i \in
\Lambda ^{n,0}(G)$, we show next that $\eta$ is closed, hence
holomorphic. Indeed,
\begin{eqnarray*} d \eta &=& \sum _i (-1)^{i+1} \omega_1 \wedge \dots
\wedge \omega_{i-1}\wedge  d \omega _i \wedge \omega_{i+1}\wedge \dots \wedge \omega_n \\
&=& \sum _i (-1)^{i+1} \omega_1 \wedge \dots \wedge \omega_{i-1}
\wedge \left(\sum_{j<i} \eta _j ^i \wedge \omega _j \right)\wedge
\omega _{i+1} \wedge \dots \wedge \omega_n =0.
\end{eqnarray*}
Since $d =\partial + \overline{\partial}\,$  and $\,
\partial \left(\Lambda^{n,0}(G)\right) \subset \Lambda^{n+1,0}(G) =0,$ it follows that
$\overline{\partial}\eta =0$, hence holomorphic.

Finally, the fact that the lattice $\Gamma $ acts on the left
implies that left invariant vector fields and $1$-forms on $G$
induce global bases of $TN$ and $T^*N$ \cite{TO}. Moreover, the
canonical projection $\pi : G \to \Gamma \backslash G$ is
holomorphic, hence the last assertion follows.
\endproof

\hfill

Another proof of triviality of the canonical bundle 
of a nilmanifold is found in Theorem 3.1 of \cite{_CG:generalized_on_nilm_}.

\hfill

On a compact K\"ahler manifold, topological triviality of the
canonical bundle implies that it is trivial holomorphically on some
finite, unramified covering of $M$. This follows  from Calabi-Yau
theorem. Indeed, by Calabi-Yau theorem, $M$ admits a Ricci-flat
K\"ahler metric (\cite{_Yau:Calabi-Yau_}). From Berger's list of
irreducible holonomies, de Rham theorem, and Cheeger-Gromoll theorem
on fundamental group of Ricci-flat manifolds, we obtain that a
finite unramified covering $\tilde M$ of $M$ is a product of compact
tori, hyperk\"ahler manifolds and simply connected Calabi-Yau
manifolds (see \cite{_Besse:Einst_Manifo_} for a detailed argument).
Therefore, $\tilde M$ has trivial canonical bundle.

On a non-K\"ahler manifold, this is no longer true. However, the
above theorem implies that the canonical bundle is holomorphically
trivial for every nilmanifold, which is never K\"ahler unless it is
a torus (see \cite{_BG_}).

For hypercomplex manifolds, $K(M,\mathcal{I})$
is always topologically trivial, which is easy to see by taking
a non-degenerate $(2,0)$-form associated with some
quaternionic Hermitian structure (Subsection
\ref{_HKT_Subsection_}). The top exterior power 
of this $(2,0)$-form trivializes  $K(M,\mathcal{I})$. However,
$K(M,\mathcal{I})$ is quite often non-trivial as
a holomorphic line bundle.

It is possible to show that $K(M,\mathcal{I})$ is non-trivial for
all hypercomplex manifolds $(M, \mathcal{I, J, K})$ such that
$(M,\mathcal{I})$ is a principal toric fibration over a base which
is a Fano manifold or orbifold
(has ample anticanonical bundle). These include
the quasiregular locally conformally hyperk\"ahler manifolds (see
\cite{_Ornea:LCHK_}), which are elliptically fibered over a contact
Fano orbifold, and compact Lie groups with the hypercomplex
structure constructed by D. Joyce (\cite{_Joyce_}), which are
torically fibered over a homogeneous rational manifold
(\cite{_Verbitsky:toric_fi_}).

Let $M\stackrel\pi \arrow B$ be such a fibration. The adjunction
formula gives $K(M)\cong \pi^* K(B)$, because the canonical bundle
of a torus is trivial. However, $\pi^* K(B)^{-N}$ has sections,
because $K(B)^{-1}$ is ample. Therefore, $K(M)$ can never be
trivial.


\section{Hypercomplex nilmanifolds and holonomy}
\label{hypercomplex_holonomy}


A manifold $(M, {\mathcal I}, {\mathcal J}, {\mathcal K})$ is called
{\bf hypercomplex} if ${\mathcal I}, {\mathcal J}, {\mathcal K}$
define integrable anticommuting complex structures on $M$ such that
${\mathcal I} {\mathcal J}= {\mathcal K}$. The operators ${\mathcal
I}, {\mathcal J}, {\mathcal K}$ define an action of the quaternion
algebra $\Bbb H$ on the tangent bundle of $M$. As Obata proved
(\cite{_Obata_}), the integrability condition of
 ${\mathcal I}, {\mathcal J}, {\mathcal K}$ is satisfied if
and only if $M$ admits a torsion-free connection $\nabla$ preserving
the quaternionic action:
\[
\nabla {\mathcal I} = \nabla {\mathcal J} = \nabla {\mathcal K} =0.
\]
Such a connection, which is necessarily unique (\cite{_Obata_}), is
called the {\bf  Obata connection} on $(M, {\mathcal I}, {\mathcal
J}, {\mathcal K})$. Setting ${\mathcal J}_1={\mathcal I}, \;
{\mathcal J}_2={\mathcal J}, \;{\mathcal J}_3={\mathcal K}$, the
Obata connection $\nabla$ is given by (see \cite{AM}):
\begin{equation} \label{Obata}
\begin{split}   \nabla_X(Y)&=
\frac{1}{12}\sum_{\alpha,\beta,\gamma} {\cal J}_{\alpha}([ {\cal
J}_{\beta}X, {\cal J}_{\gamma}Y ]+[{\cal J}_{\beta}Y,{\cal
J}_{\gamma}X ])
\\ & + \frac{1}{6} \sum_{\alpha=1}^3{\cal J}_{\alpha}([{\cal J}_{\alpha}X, Y
]+[{\cal J}_{\alpha}Y, X ])+ \frac{1}{2}[X,Y] ,
\end{split}
\end{equation}
$X,Y \in \frak X(M)$, where ${\alpha},{\beta},{\gamma}$ is a cyclic
permutation of $1,2,3$.

%
We consider next hypercomplex nilmanifolds.

\hfill

\definition
A hypercomplex structure ${\mathcal I}, {\mathcal J}, {\mathcal K}$
on a Lie group $G$ is called {\bf left-invariant} when left
translations are holomorphic with respect to the complex structures
${\mathcal I}, {\mathcal J}$ and $ {\mathcal K}$. Let $N=\Gamma
\backslash G$ be a nilmanifold, with $G$ a Lie group equipped with a
left-invariant hypercomplex structure.  The quotient $N=\Gamma
\backslash G $ inherits a hypercomplex structure. In such situation,
we say that $N$ is a {\bf hypercomplex nilmanifold}.

\hfill

Let $\Hol(\nabla)$ be the holonomy group associated with the Obata
connection~$\nabla$. Since $\nabla$ preserves the quaternionic
structure, $\Hol(\nabla)\subset GL(n, {\Bbb H})$. We define the
determinant of $h\in GL(n, {\Bbb H})$ in the following way. Let
$V\cong {\Bbb H}^n$ be the vector space over $\Bbb H$, and $V_{
I}^{1,0}$ the same space considered as a complex space with the
complex structure $I$ induced by ${\mathcal I}$. The Hodge
decomposition gives $V \otimes _{\R}\C \cong V_{ I}^{1,0}\oplus V_{
I}^{0,1}$. The top exterior power $\Lambda^{2n,0}_{
I}(V):=\Lambda^{2n}(V_{ I}^{1,0})\cong \C$ is equipped with a
natural real structure:
\begin{equation}\label{_real_stru_from_J_Equation_}
\eta \arrow { J}(\bar\eta)
\end{equation}
for $\eta \in \Lambda^{2n,0}_{ I}(V)$ (the quaternions ${ I}$ and
${J}$ anticommute, hence ${ J}$ exchanges $\Lambda^{p,q}_{ I}(V)$
with $\Lambda^{q,p}_{ I}(V)$). Since the real structure on
$\Lambda^{2n,0}_{ I}(V)$ is constructed from the quaternion action,
any $h\in GL(V,{\Bbb H})$ preserves this real structure. Let
$\det(h)$ denote the action induced by $h$ on $\Lambda^{2n,0}_{
I}(V)\cong \C$. Then $\det (h)\in \R$, as the above argument
implies. This defines a homomorphism
\[
\det:\; GL(n, {\Bbb H}) \arrow \R^*
\]
to the multiplicative group of non-zero real numbers, which is
clearly positive since $GL(n, {\Bbb H})$ is connected. Let $SL(n,
\H)\subset GL(n, {\Bbb H})$ be the kernel of $\det$, $K(M,{\mathcal
I})$  the canonical bundle of $(M,{\mathcal I})$ and $\nabla_K$  the
connection on $K(M,{\mathcal I})$ induced by the Obata connection
$\nabla$. Given $h\in \Hol(\nabla)$, the corresponding
transformation in $\Hol(\nabla _K)$ acts on sections of
$K(M,{\mathcal I})$ by multiplication by $\det(h)$, hence
\[\Hol(\nabla _K)= \{ \det(h) : h \in \Hol (\nabla)\}. \]
Therefore, $\nabla_K$ has trivial holonomy if and only if
$\Hol(\nabla)\subset SL(n, {\Bbb H})$.
Moreover, the last condition implies that $K(M,{\mathcal I})$ is
holomorphically trivial (see \cite{_Verbitsky:canoni_}, Claim 1.2).
We show in \ref{_SL(_n_H_)_Corollary_} that the converse of this
statement holds in the case that  $M$ is a nilmanifold, thereby
giving an affirmative answer to a question raised in
\cite{_Verbitsky:canoni_}.  The proof of this corollary makes use of
\ref{_canoni_trivial_Theorem_} and the next result:

\hfill


\theorem\label{_canon_tri_then_SL_Theorem_} Let $(N, {\mathcal I},
{\mathcal J}, {\mathcal K})$ be a hypercomplex nilmanifold, $\dim_\C
N=2n$, and $\eta$  a holomorphic, left-invariant section of the
canonical bundle $ \Lambda^{2n,0}(N,{\mathcal I})$. Then
$\nabla\eta=0$, where $\nabla$ is the Obata connection.

\hfill

{\bf Proof:} Since the Obata connection is torsion-free, $d\eta =
\Alt(\nabla \eta)$, where $\Alt= \bigwedge:\; \Lambda^{2n}(M)
\otimes \Lambda^1(M)\arrow\Lambda^{2n+1}(M)$ denotes the exterior
product. Since $\eta$ is holomorphic, $\bar\6\eta=0$. The map $\Alt$
restricted to $\Lambda^{2n,0}(M) \otimes \Lambda^{0,1}(M)$ is an
isomorphism; therefore,
 $\nabla^{0,1}\eta=0$.

Any left-invariant section of $\Lambda^{2n,0}_I(N)$ is holomorphic,
because such a section is unique, up to a constant multiplier.
Therefore, $J(\bar\eta)\in \Lambda^{2n,0}_I(N)$ is holomorphic. This
gives
\begin{equation} \label{_nabla_of_J(eta)_Equation_}
\nabla^{0,1}J(\bar\eta)=0.
\end{equation}
Since $\nabla$ commutes with $J$, \eqref{_nabla_of_J(eta)_Equation_}
implies $\nabla^{0,1}\bar\eta=0$. However,
$\nabla^{0,1}\bar\eta=\overline{\nabla^{1,0}\eta}$, and this gives
$\nabla^{1,0}\eta=0$. We proved that $\nabla^{0,1}\eta +
\nabla^{1,0}\eta = \nabla\eta =0$.
\endproof

 \endproof

\hfill

Comparing \ref{_canon_tri_then_SL_Theorem_} with
\ref{_canoni_trivial_Theorem_}, we obtain the following important
corollary:

\hfill

\corollary\label{_SL(_n_H_)_Corollary_} Let $(N, {\mathcal I},
{\mathcal J}, {\mathcal K})$ be a hypercomplex nilmanifold. Then
$\Hol(\nabla)\subset SL(n, {\Bbb H})$, where $\Hol(\nabla)$ is the
holonomy of the Obata connection.

\hfill

\noindent {\bf Proof:} \ref{_canoni_trivial_Theorem_} implies that
$\Lambda^{2n,0}(N,{\mathcal I})$ has a holomorphic section and by
\ref{_canon_tri_then_SL_Theorem_}, $\Hol(\nabla)\subset SL(n, {\Bbb
H})$ where $\nabla$ is the Obata connection.
\endproof

\hfill

As a consequence of the above result it follows that the Obata
connection on any hypercomplex nilmanifold is Ricci flat.

\hfill

\corollary \label{ricciflat} Let $(N, {\mathcal I}, {\mathcal J},
{\mathcal K})$ be a hypercomplex nilmanifold. Then the Ricci tensor
of the Obata connection vanishes.

\hfill

We give two proofs of this corollary; the first one is a consequence
of Theorem 5.6 in \cite{AM} and the second one makes use of Lemma
3.2 in \cite{_Bar_}. In both proofs  $\nabla$ denotes the Obata
connection of the left-invariant hypercomplex structure on $G$,
where  $N=\Gamma \backslash G$.

\hfill

\noindent
 {\bf  First Proof:}    It follows from \ref{_SL(_n_H_)_Corollary_} that
$\Hol(\nabla)\subset SL(n, {\Bbb H})$. It was proved in \cite{AM},
Theorem 5.6\footnote{Theorem 5.6 in \cite{AM} holds for  $n>1$. For
$n=1$ it still holds if we  assume that $W_+=0$, where $W_+$ is the
self-dual part of the Weyl tensor $W$. This assumption is immediate
for hypercomplex manifolds, because the hypercomplex structure gives
a parallel trivialization of the bundle $\Lambda^+(M)$.}, that for a
simply connected
 hypercomplex manifold of dimension $4k$, $k>1$, the Obata connection $\nabla$ satisfies $\Hol(\nabla)\subset SL(n, {\Bbb
H})$ if and only if the Ricci tensor of $\nabla$ vanishes.
Therefore,   the Ricci tensor of $\nabla$ vanishes on $G$, hence it
vanishes on $N$.
\endproof

\hfill

\noindent
 {\bf  Second Proof:}
 Let $\frak g$ be the Lie algebra of $G$ and set ${\mathcal J}_1={\mathcal
I}, \; {\mathcal J}_2={\mathcal J}, \;{\mathcal J}_3={\mathcal K}$.
According to  Lemma 3.2 in \cite{_Bar_},
\[  Ric \equiv 0   \quad \text{ if and only if } \quad
 \text{tr}\left(\nabla_{[X_1,X_2]} \right) = 0,  \forall \, X_1, X_2 \in \frak g  . \]
 The first step is to show that:
\begin{equation}  \label{1st} \text{tr}\left(\nabla_{[X_1,X_2]} \right) =
\text{tr}\left({\cal J}_{\alpha}\, \text{ad}_{{\cal
J}_{\alpha}[X_1,X_2]}\right), \quad X_1, X_2 \in \frak g, \; \alpha
=1,2,3.\end{equation}

 We compute the trace of
$\nabla_{[X_1,X_2]}$ (recall \eqref{Obata}):
\begin{equation} \label{traza}
\begin{split}
& \text{tr}\left(\nabla_{[X_1,X_2]} \right) = \frac{1}{6} \,
\text{tr}\left(\sum_{\alpha =1}^3{\cal J}_{\alpha}\,
\text{ad}_{{\cal
J}_{\alpha}[X_1,X_2]}\right)\\
&+\frac{1}{12}\, \text{tr} \left( \sum_{\alpha,\beta,\gamma} \left(
{\cal J}_{\alpha}\text{ad}_ {{\cal J}_{\beta}[X_1,X_2]} {\cal
J}_{\gamma}-{\cal J}_{\alpha}ad_{{\cal J}_{\gamma}[X_1,X-2]}{\cal
J}_{\beta}\right)\right) ,
\end{split}
\end{equation}
where $\alpha, \beta, \gamma$ is a cyclic permutation of $1,2,3$.

Since $$\text{tr}\, \left({\cal J}_{\alpha}\text{ad}_ {{\cal
J}_{\beta}[X_1,X_2]} {\cal J}_{\gamma}\right)=\text{tr}\left({\cal
J}_{\beta}\text{ad}_ {{\cal J}_{\beta}[X_1,X_2]}\right)$$ and
$$\text{tr}\left({\cal J}_{\alpha}\text{ad}_ {{\cal
J}_{\gamma}[X_1,X_2]} {\cal J}_{\beta}\right)= -\text{tr}\left({\cal
J}_{\gamma}\text{ad}_ {{\cal J}_{\gamma}[X_1,X_2]}\right)$$ it
follows that:
$$\text{tr}\left(\nabla_{[X_1,X_2]}\right)= \frac{1}{3}\, \text{tr}\left(\sum_{\alpha=1}^3{\cal J}_{\alpha}\, \text{ad}_{{\cal
J}_{\alpha}[X_1,X_2]}\right).$$

We show next that:
\begin{equation}\label{trJad_J} \text{tr}\left({\cal
J}_{\alpha}\,\text{ad}_{{\cal J}_{\alpha}[X_1,X_2]}\right) \; \text{
is independent of  } \; \alpha=1,2,3.\end{equation} Set
$X=[X_1,X_2]$ and let $Y \in \frak g$. Observe that:
\begin{equation} \label{trace}\begin{split}
\text{tr}\left({\cal J}_{\alpha}\, \text{ad}_{{\cal J}_{\alpha}X
}\right)&= \text{tr}\left( \text{ad}_{{\cal J}_{\alpha}X } {\cal
J}_{\alpha}\right)=-\text{tr}\left({\cal J}_{\gamma}\,
\text{ad}_{{\cal J}_{\alpha}X } {\cal J}_{\alpha} {\cal
J}_{\gamma}\right)\\&= \text{tr}\left({\cal J}_{\gamma}\,
\text{ad}_{{\cal J}_{\alpha}X } {\cal J}_{\beta}\right).\end{split}
\end{equation} The integrability of ${\cal J}_{\gamma}$
gives:
$$ {\cal J}_{\gamma}[{\cal J}_{\alpha}X, {\cal J}_{\beta}Y]=[{\cal J}_{\beta}X,{\cal J}_{\beta}Y]- [{\cal J}_{\alpha}X, {\cal J}_{\alpha}Y]-
{\cal J}_{\gamma}[{\cal J}_{\beta}X, {\cal J}_{\alpha}Y] ,$$ which
implies that: \begin{equation*} \text{tr}\left({\cal J}_{\gamma}\,
\text{ad}_{{\cal J}_{\alpha}X } {\cal J}_{\beta}\right)=
\text{tr}\left( \text{ad}_{{\cal J}_{\beta} X } {\cal
J}_{\beta}\right) - \text{tr}\left( \text{ad}_{{\cal J}_{\alpha}X }
{\cal J}_{\alpha}\right)-  \text{tr}\left({\cal J}_{\gamma}\,
\text{ad}_{{\cal J}_{\beta}X } {\cal J}_{\alpha}\right).
\end{equation*}
Using \eqref{trace} we obtain:
$$ \text{tr}\left( {\cal
J}_{\alpha} \text{ad}_{{\cal J}_{\alpha}X } \right)= \text{tr}\left(
\text{ad}_{{\cal J}_{\beta} X } {\cal J}_{\beta}\right) -
\text{tr}\left( \text{ad}_{{\cal J}_{\alpha}X } {\cal
J}_{\alpha}\right)+\text{tr}\left( {\cal J}_{\beta} \text{ad}_{{\cal
J}_{\beta}X } \right),$$ or equivalently,
$$ 2\, \text{tr}\left( {\cal
J}_{\alpha} \text{ad}_{{\cal J}_{\alpha}X } \right)= 2\,
\text{tr}\left( {\cal J}_{\beta} \text{ad}_{{\cal J}_{\beta}X }
\right),$$ and \eqref{trJad_J} follows. This implies \eqref{1st},
which together
 with \ref{ric0} imply that the
corollary holds on $G$, hence on $N=\Gamma \backslash G$.
\endproof

\hfill

\remark Notice that the converse of \ref{_SL(_n_H_)_Corollary_} is
not necessarily true. Indeed, the vanishing of the Ricci curvature
is equivalent to the flatness of the canonical bundle $K(N)$ of $N$.
However, it might  have global monodromy, as it happens in the case
of the Hopf surface. \ref{_SL(_n_H_)_Corollary_} implies that (for a
nilmanifold), $K(N)$ is trivial, both locally and globally.


\section{Quaternionic Hermitian structures on nilmanifolds}


\subsection{HKT structures on abelian nilmanifolds}
\label{_HKT_Subsection_}

Let $(M, {\mathcal I}, {\mathcal J}, {\mathcal K})$  be a
hypercomplex manifold. A quaternionic Hermitian metric $g$ on $M$ is
a Riemannian metric which is Hermitian with respect to ${\mathcal
I}, {\mathcal J}$ and $ {\mathcal K}$. This is equivalent to $g$
being $SU(2)$-invariant with respect to the $SU(2)$-action generated
by the group of $SU(2)\cong SU({\Bbb H}, 1)$ of unitary quaternions,
\[ SU({\Bbb H}, 1)= \{ a+ bI + cJ + dK\ | \  a^2+b^2+c^2+d^2=1\}.
\]
Starting from an arbitrary Riemannian metric and averaging over
$SU(2)$, we obtain a quaternionic Hermitian metric. Therefore, such
metric always exists.

\hfill

\definition
Let $(N, {\mathcal I}, {\mathcal J}, {\mathcal K})$  be a
hypercomplex nilmanifold, where $N=\Gamma \backslash G$. When the
Lie group
 $G$ is equipped with a left-invariant quaternionic Hermitian metric
 we say that the hypercomplex nilmanifold $N$
 with the induced
 metric is a {\bf quaternionic Hermitian nilmanifold}.

\hfill

Let $(M, {\mathcal I}, {\mathcal J}, {\mathcal K})$ be a
hypercomplex manifold, $g$ a quaternionic Hermitian metric, and
$\Omega$ a 2-form on $M$ constructed from $g$ as follows:
\begin{equation}\label{_Omega_Equation_}
\Omega := g( \cdot, {\mathcal J}\cdot) +\1 g( \cdot, {\mathcal
K}\cdot )
\end{equation}
Then, $\Omega$ is a $(2,0)$-form on $(M,{\mathcal I})$ as an
elementary linear-algebraic argument implies
(\cite{_Besse:Einst_Manifo_}).

 The hyperk\"ahler
condition can be written down as $d\Omega=0$
(\cite{_Besse:Einst_Manifo_}). The HKT condition is weaker:

\hfill

\definition\label{_HKT_Definition_}
A quaternionic Hermitian metric is called an HKT-metric if
\begin{equation}\label{_HKT_intro_Equation_}
\6(\Omega)=0,
\end{equation}
where $\6:\; \Lambda^{2,0}_{\mathcal I}(M) \arrow
\Lambda^{3,0}_{\mathcal I}(M)$ is the Dolbeault differential on $(M,
{\mathcal I})$, and $\Omega$ the $(2,0)$-form on $(M, {\mathcal I})$
constructed from $g$ as in \eqref{_Omega_Equation_}.

\hfill

\definition
Let $(N, {\mathcal I}, {\mathcal J}, {\mathcal K})$  be a
hypercomplex nilmanifold, where $N=\Gamma \backslash G$. When the
Lie group
 $G$ is equipped with a left-invariant HKT-metric we say that the hypercomplex nilmanifold $N $
 with the induced
 metric is an {\bf HKT nilmanifold}.

\hfill

\remark It has been shown in \cite{_Fino_Gra_} that existence of
any HKT-metric on $(N, {\mathcal I}, {\mathcal J}, {\mathcal K})$,
compatible with a left-invariant hypercomplex structure implies
existence of a left-invariant one.

\hfill

\definition
 A hypercomplex nilmanifold $(N,{\mathcal I}, {\mathcal J},
{\mathcal K})$ with $N=\Gamma \backslash G$ is called {\bf abelian}
when ${\mathcal I}, {\mathcal J}, {\mathcal K}$ are induced by
left-invariant abelian complex structures  on $G$.

\hfill

In \cite{_Dotti_Fino:HKT_}, it was shown that for each invariant
abelian hypercomplex structure on a Lie group, any left-invariant
quaternionic Hermitian metric is HKT. This implies that any abelian
hypercomplex nilmanifold is HKT. We show next that as a consequence
of \ref{_canoni_trivial_Theorem_}  the converse of this result
holds. The case of $2$-step nilmanifolds was proved in
\cite{_Dotti_Fino:HKT_}.

\hfill

\theorem \label{_HKT_implies_abelian_} Let $(N,{\mathcal I},
{\mathcal J}, {\mathcal K},g)$ be a nilmanifold admitting an
HKT-structure. Then it is abelian.

\hfill

The proof of the above theorem will follow from a Hard Lefschetz
isomorphism on the Dolbeault cohomology of $(N, {\mathcal I})$.

\hfill

\proposition  \label{_nil_hard_lefsch_} Let $(N,{\mathcal I},
{\mathcal J}, {\mathcal K}, g)$ be an HKT nilmanifold and $\Omega$
the corresponding $(2,0)$-form with respect to $\mathcal I$ (see
\eqref{_Omega_Equation_}). Then,
$$ L_{\Omega}^{n-i}: H^{i,0}_{\partial}(N, {\mathcal I})\rightarrow
H^{2n-i,0}_{\partial}(N, {\mathcal I})$$ is an isomorphism, where
$L_{\Omega}\left([\gamma]\right)= [{\Omega}\wedge \gamma]$.

\hfill

\noindent {\bf Proof:} Let ${\omega_1}, \dots, {\omega _{2n}}$ be a
basis of invariant $(1,0)$-forms on $N$ as in the proof of
\ref{_canoni_trivial_Theorem_}. Then $
\bar{\eta}=\bar{\omega}_1\wedge \dots \wedge \bar{\omega} _{2n}$ is
an invariant section of the line bundle $\Lambda^{0,2n}(N,{\mathcal
I})$. Therefore, $\bar{\Omega}^n$ is proportional (with a constant
factor) to $\bar{\eta}$.  Let $\theta$ be the $(1,0)$-form defined
by
$$ \partial \bar{\Omega}^n = \theta \wedge  \bar{\Omega}^n .    $$
Since $\bar{\eta}$ is closed, d$\,\bar{\Omega}^n =0$, hence
$\partial \bar{\Omega}^n =0$ and it follows that $\theta =0$. This
says that the Dolbeault complex of the square root of the canonical
bundle $K(N, \mathcal I)$ determined by the trivialization induced
by $\Omega ^n$ is identified with the complex $\left(
\Lambda^{*,0}(N, \mathcal I), \partial \right) $. The proposition
now follows from Theorem~10.2 in~\cite{_Verbitsky:HKT_}.
\endproof

\hfill

 As a consequence of the above result and Lemma 9 in
\cite{_Console_Fino_}  we obtain:

\hfill

\corollary \label{_alg_hard_lefsch_} Let $(N,{\mathcal I}, {\mathcal
J}, {\mathcal K}, g)$ be an HKT nilmanifold, with $N=\Gamma
\backslash G$. Then,
$$ L_{\Omega}^{n-i}: H^{i,0}_{\partial}(\frak g _{\Bbb C}, { I})\rightarrow
H^{2n-i,0}_{\partial}(\frak g _{\Bbb C}, { I})$$ is an isomorphism,
where $\frak g _{\Bbb C}$ is the complexification of the  Lie
algebra of $G$.

\hfill

\noindent {\bf Proof of \ref{_HKT_implies_abelian_}:} The aim is to
show that $\frak g ^{1,0}$ is abelian. If $\frak g ^{1,0}$ were not
abelian, an analogous argument to that in \cite{_BG_} would give
that
$$ L_{\Omega}^{n-1}: H^{1,0}_{\partial}(N, {\mathcal I})\rightarrow
H^{2n-1,0}_{\partial}(N, {\mathcal I}) $$ is not surjective; this
contradicts \ref{_alg_hard_lefsch_}. Therefore, $\frak g ^{1,0}$
must be abelian. Repeating the argument with $ {\mathcal J}$ and $
{\mathcal K}$  the theorem follows.
\endproof

\hfill

The next corollary is a straightforward consequence of
\ref{_HKT_implies_abelian_} (compare with Theorem 3.1 in
\cite{_Dotti_Fino:HKT_}):

\hfill

\corollary The hypercomplex structure of a left-invariant HKT-metric
on  a   nilpotent Lie group admitting a lattice is abelian.

\subsection{Quaternionic  balanced metrics}

\definition A quaternionic Hermitian metric $g$ on a hypercomplex
manifold  is called {\bf quaternionic balanced} if it is balanced
with respect to all complex structures.

\hfill

\proposition \label{quat_balanced} Let $(N, {\mathcal I},
{\mathcal J}, {\mathcal K}, g)$ be a quaternionic Hermitian
nilmanifold such that the hypercomplex structure is abelian. Then
$g$ is quaternionic balanced.

\hfill

 \noindent {\bf  Proof:} Let $N=\Gamma \backslash G$; we still
 denote by
${\mathcal I}, {\mathcal J}, {\mathcal K}, g$ the induced
left-invariant quaternionic Hermitian  structure on $G$. As shown
in \cite{_Dotti_Fino:HKT_}, $g$ is HKT. Therefore, the Bismut
connections associated with ${\mathcal I}, {\mathcal J},
 {\mathcal K}$ are equal (this is one of the alternative
definitions of HKT-structures, see  \cite{_Gra_Poon_} for
details). Denote the Bismut connection of $N$ by $\nabla^B$.

 Since $G$ is nilpotent,  formula \eqref{abel_Lee} implies that the Lee
 form $\theta _J$ corresponding to $(J,g)$ is given by
 \begin{equation}  \theta_J (X)= \text{tr}  \left(  \frac 12 \, J \nabla^B_{JX}
 \right), \qquad X \in \frak g,
 \end{equation}
where $J$ is the  complex structure on $\frak g$ induced by
$\mathcal J$. We show next that $\text{tr}  \left(   \, J
\nabla^B_{JX}\right)=0$. Let $X_1, IX_1, JX_1, K X_1 ,\dots , X_n,
IX_n, JX_n, KX_n$ be an orthonormal basis of $\frak g$. From the
definition of the Bismut connection it follows
\[\nabla^BI=\nabla^BJ=\nabla^BK=0.\]
Then
\begin{equation*}\begin{split}
\text{tr}\left(J \nabla^B_{JX}\right)=&\sum_{j=1}^n g\left(J
\nabla^B_{JX}X_j, X_j\right) + \sum_{j=1}^n g\left(J
\nabla^B_{JX}IX_j, IX_j\right) \\ & +\sum_{j=1}^n g\left(J
\nabla^B_{JX}JX_j, JX_j\right)+\sum_{j=1}^n g\left(J
\nabla^B_{JX}KX_j, KX_j\right)\\
=&\sum_{j=1}^n g\left(J \nabla^B_{JX}X_j, X_j\right) +
\sum_{j=1}^n g\left(JI \nabla^B_{JX}X_j, IX_j\right) \\ &
+\sum_{j=1}^n g\left( \nabla^B_{JX}JX_j, X_j\right)+\sum_{j=1}^n
g\left(JK
\nabla^B_{JX}X_j, KX_j\right)\\
=&\sum_{j=1}^n g\left(J \nabla^B_{JX}X_j, X_j\right) -\sum_{j=1}^n
g\left(IJ \nabla^B_{JX}X_j, IX_j\right) \\ & +\sum_{j=1}^n g\left(
\nabla^B_{JX}JX_j, X_j\right)-\sum_{j=1}^n g\left(KJ
\nabla^B_{JX}X_j, KX_j\right)\\
=&\sum_{j=1}^n g\left( \nabla^B_{JX}JX_j, X_j\right) -\sum_{j=1}^n
g\left(\nabla^B_{JX}J X_j,  X_j\right) \\ & \sum_{j=1}^n g\left(
\nabla^B_{JX}JX_j, X_j\right)-\sum_{j=1}^n g\left(
\nabla^B_{JX}JX_j, X_j\right) =0,
\end{split}\end{equation*}
therefore, $\theta_J(X)=0$ and $g$ is balanced with respect to $J$.
The same proof holds for $I$ and $K$.
\endproof

\hfill

As a consequence of \ref{_HKT_implies_abelian_} and
\ref{quat_balanced} we obtain:

\hfill

\begin{corollary}  Let $(N,{\mathcal I},
{\mathcal J}, {\mathcal K},g)$ be an
 HKT nilmanifold. Then $g$ is quaternionic balanced.
\end{corollary}

\subsection{A family of non-abelian hypercomplex
 nilmanifolds}
\label{examples}

We end this section by exhibiting a family of hypercomplex
nilmanifolds which do not admit HKT metrics. This will follow from
\ref{_HKT_implies_abelian_} since such hypercomplex nilmanifolds are
not abelian.

Let $A$ be
 a finite dimensional  associative  algebra  and $\frak a \frak f \frak f (A)$ the Lie algebra $A \oplus A$
with Lie  bracket  given as follows:
\[   [(a,b),(a',b')]=(aa'-a'a,ab'-a'b), \hspace{1cm} a,b,a',b' \in A . \]
These Lie algebras have been considered in \cite{BD2}. We note that
$\frak a \frak f \frak f (A)$ is a nilpotent Lie algebra if and only
if $A$ is nilpotent as an associative algebra.

 Let $J$ be the endomorphism  of $\frak a \frak f
\frak f (A)$ defined by
\begin{equation} J(a,b)=(b,-a),    \hspace{1cm} a,b  \in A .
\label{jaff1} \end{equation}
 A computation shows  that $J$ defines a complex structure
on $\frak a \frak f \frak f (A)$.   Furthermore, if
 one assumes the algebra $A$ to be a complex associative  algebra,
this extra assumption allows us to equip $\frak a \frak f \frak f
(A)$ with a pair of anti-commuting
 complex structures.  Indeed,
 the endomorphism $K$ on  $\frak a \frak f \frak f (A)$
 defined by $K(
a,b)=(-ia, ib)$ for $a,b\in A$ is a complex structure anticommuting
with $J$, hence, $J$ and $K$ define a hypercomplex structure by
setting $I=JK$. Moreover, the hypercomplex structure is abelian if
and only if $A$ is  commutative. It then follows that the simply
connected Lie groups with Lie algebra $\frak a \frak f \frak f (A)$,
where $A$ is a complex associative non-commutative algebra carry
non-abelian hypercomplex structures. In particular, let $A_k$ be the
algebra of $k\times k$ strictly upper triangular matrices with
complex entries and Aff$(A_k)$ the simply connected Lie group with
Lie algebra $\frak a \frak f \frak f (A_k)$, which is $(k-1)$-step
nilpotent. Since the  structure constants with respect to the
standard basis of $\frak a \frak f \frak f (A_k)$ are integers,
there exists a lattice   $\Gamma _k$ in Aff$(A_k)$, thus the
hypercomplex nilmanifold $\, N_k=\Gamma_k \backslash$Aff$(A_k)\, $
does not carry an HKT-metric.

\hfill

\noindent{\bf Acknowledgements:} Misha Verbitsky thanks Geo
Grantcharov and Ma\-xim Kontsevich for interesting discussions of
Bogomolov-Tian-Todorov theorem. We are grateful to Maxim Kontsevich
for the reference to \cite{_Ghys_}.

\hfill

\hfill

\hfill

\small{

\noindent{\sc Mar\'\i a L. Barberis, \ \
Isabel G. Dotti\\
FaMAF - Universidad Nacional de C\'ordoba, CIEM - CONICET, \\
Ciudad Universitaria, 5000 C\'ordoba, Argentina}\\
\tt barberis@mate.uncor.edu, \  \
 idotti@mate.uncor.edu\\[4mm]
\noindent {\sc Misha Verbitsky\\
{\sc  Institute of Theoretical and
Experimental Physics \\
B. Cheremushkinskaya, 25, Moscow, 117259, Russia }\\
\tt   verbit@mccme.ru }

}

\end{document}